\documentclass[aop]{imsart}
\usepackage{graphicx}                      

\usepackage{amsfonts}
\usepackage[margin=1.0 in]{geometry}

\RequirePackage{amsthm,amsmath}  

\usepackage{enumitem}
\usepackage{amssymb}

\startlocaldefs
\numberwithin{equation}{section}
\theoremstyle{plain}
\newtheorem{thm}{Theorem}[section]

\endlocaldefs

\begin{document}

\begin{frontmatter}
\title{
         Chernoff bounds for branching random walks
}

\begin{aug}
\affiliation{Changqing Liu  \\ ClinTFL Ltd. \\  c.liu@ClinTFL.com}
\end{aug}

\begin{abstract}
Concentration inequalities, which have proved very useful in a variety of fields, provide fairly tight bounds on large deviation probabilities while central limit theorem (CLT) describes the asymptotic distribution around the mean (at the  $\sqrt{n}$ scale).
Harris (1963) conjectured that for a supercritical branching random walk (BRW) of i.i.d offspring and i.i.d displacements, positions of individuals in $nth$ generation approach to Gaussian distribution --- central limit theorem. This conjecture was later proved by Stam  (1966) and Kaplan \& Asmussen (1976). Refinements and extensions followed. However, to the best of our knowledge, there is no corresponding existing work on concentration inequalities for BRWs. In this note, we propose a new definition of BRW, providing a more general framework. Owing to this definition, a Chernoff-type (subgaussian) bound for BRWs follows directly from the Chernoff bound for random walk. The relation between RW (random walk) and BRW is discussed.
\end{abstract}


\end{frontmatter}

\section{\textbf Introduction}

\subsection {Concentration inequalities}
Let $(S_i)_{i=1, 2, ...}$ be a random process of one-dimensional random walk on the real line. $S_n$ can be formulated as $S_n=X_1 + X_2 + \cdots + X_n$; $X_i=S_i - S_{i-1}$. Under certain conditions, among which independence is crucial, $S_n$ converges to a Gaussian distribution in a $\sqrt{n}$ neighborhood of its mean --- central limit theorem (CLT). For deviation beyond this neighborhood, i.e., when $|S_n - \bold E (S_n)| \gg \sqrt{n}$, Chernoff inequality (\cite{Chernoff 1952} \cite{Hoeffding 1963}),  also referred to as Chernoff bound,
\begin{equation}\label{chernoff bound}
         \Pr(|S_n - \bold E(S_n)| \geq \lambda) \leq e^{- c\frac{\lambda^2}{n}}
\end{equation}
describes how unlikely this occurrence is. In this note, by ``Chernoff bound'' we refer to the tail probability inequalities of the form (\ref{chernoff bound}) -- the subgaussian (square-exponential) deviation bounds. Chernoff inequalities are extended to the setting of bounded martingale difference sequences (e.g., $E(X_n \ | \ S_{n-1}) = 0$), and are then referred to as Azuma-Hoeffding inequality  \cite{Azuma}  (see \cite{Chung}, \cite{McDiarmid89} for the survey and references therein). Refinements and extensions followed (e.g. \cite{Talagrand 1995}, \cite{McDiarmid89}, \cite{Bentkus 2004}, \cite{bentkus 2008}, \cite{Fan & Grama & Q. Liu 2012}, \cite{Kontorovich 2014}). Known collectively as concentration inequalities, they prove extremely useful and have a wide variety of applications in computer science, combinatorics, information theory (see e.g. \cite{Dubashi&Panconesi 2009}, \cite{Alon2000},  \cite{Raginsky&Sason 2013}, \cite{Motwani&Raghavan 1995}).

\subsection{Branching random walk (BRW)}
In the literature, a branching random walk on the real line is described as follows (see, for example, \cite{Biggins_fst_lst} \cite{Biggins_Mtgl} \cite{Biggins_Chernoff}). In generation zero, an initial particle at the origin on the real line $\mathbf{R}$. It splits into a random number of child particles who form generation one. The children's displacements, relative to their parent, correspond to a point process on $\mathbf{R}$. The children in turn split too to form the second generation, and so on. If the average split number (branching factor) is greater than one, with positive probability the number of the descendants grows exponentially through generations. Current BRW studies typically address models where the offspring's behavior is independent of that of their previous generation (e.g. \cite{Biggins_Mtgl}, \cite{Biggins_Chernoff} \cite{BIGGINS_et1997}, \cite{Hu&Shi2009}, \cite{Gao_BRW2014}). The law of large numbers and central limit theorem type results about the distribution of position are established, under some conditions of independence (for instance, i.i.d of branching and walking). However, unlike random walk, Chernoff bound is not known so far even in the case of i.i.d aforementioned, while the minimal (and maximal) is studied by many.

Considering i.i.d. offspring (and hence independent of the parent's position), and i.i.d. displacement, Harris \cite{Harris63} conjectured that the distribution of the descendants' position of the $n$th generation approaches Gaussian distribution. This conjecture was proved by \cite{Stam66}, \cite{Kaplan_II}. Its extended generation-dependent versions, where offspring and displacement distribution are dependent on generation $n$, were proved by \cite{Biggins_CLT}, \cite{Klebaner82}, to mention a few. Problems where offspring and displacement are dependent of parents' positions are studied in adhoc approaches though (e.g. \cite{Yoshida}). Concerning the deviation from the expectation, \cite{Biggins_Chernoff} showed that for any fixed $\delta > 0$, the proportion of particles whose scaled position is less than $\mu - \delta$ (or greater than $\mu + \delta$) converges almost surely to zero. In another direction of estimating the rareness of the large deviation, extremum is well studied (e.g. \cite{DEKKING1991} \cite{MCDIARMID_95}, \cite{BACHMANN2000}, \cite{BRAMSON&Zeitouni 2009} \cite{Hu&Shi2009}, \cite{Addario2009}  \cite{Aidekon 2013} ). \cite{DEKKING1991} showed tightness for $M_n - E M_n$, where $M_n$ is the minimal position in $n$th generation. \cite{MCDIARMID_95} gave probability bound for the deviation of $M_n$, i.e. $\Pr(M_n - Med_{n} > x) < e^{-\delta x}$, where $Med_n$ is the median of $M_n$. While probability bounds for deviations at the $n$-scale (i.e., for fixed $\epsilon > 0$ under $n^{-1}$ scaling) are by now established, little is known about how sparse the population is at distances of order $\sqrt{n}$ from the mean (position), in the concentration-inequality sense familiar from random walks.

\section {Chernoff bound for BRW}
Throughout we consider BRW on $\mathbf{R}$. We introduce a different definition of BRW from a new perspective, which allows us to treat BRW in a more general way in three respects: First, in our framework both subcritical and supercritical BRW are treated without distinction; while traditionally subcritical BRW is considered trivial, since a branching random walk process almost surely ends with zero population when $n$ is large, we observed that the probability space for the survival BRW processes is well defined with infinitely large number of ancestors. Second, the underlying random walks between siblings are not assumed to be independent. Third, branching factor (birth rate) is not assumed identical across generations and siblings. The only major requirement is the independence between birth-rate and birth-place.

Starting from one initial ancestor (1st parent), a realization of branching random walk process is a random rooted tree. Each node in the tree is associated with a position which equals to the sum of the displacements of its previous parents and the displacement of itself, $u$. As a parent located at $u$, this node in turn produces some number of children with $u$ as the birth-place. Infinite BRW processes with such initial ancestors make a forest, a probability space of the BRW. In other words, the probability space of the BRW can be interpreted as a forest grows from infinite roots ---initial ancestors. Let 
\[
      S_n = X_1 + X_2 + \cdots + X_n
\]
be the position of a leaf in the $n$th generation, where $X_{i}$ is displacement (step size relative from the birth place) of its $ith$ parent. $(X_1, X_2, ..., X_n)$ is called a  {\it spine} (a path from the origin to $S_n$). The central limit theorem for BRW states that almost surely a randomly chosen tree has a Gaussian ``canopy"; namely, $S_n$ on the tree has Gaussian distribution. Let $u_i$ be an individual's position, i.e. $X_1 + \cdots+ X_i$, at generation $i$. At position $u_i$, the mean of birth rate is $m_i(u)$ (generally the mean of birth rate may be dependent on the birth place).

In this paper, we study BRW from a new perspective. Traditionally, BRW is viewed as Galton-Watson tree weighted by spatial displacements and its study relies heavily on generating-function techniques. In contrast, we take BRW as RW $+$ point process---maybe called ``RWB''. In a random walk, viewed as a special BRW, a parent produces exactly one child and is replaced by it. The intensity measure of the displacement point process $\{\xi_1\}$ is 
\[
          \bold E \left[ \mathbf{1}_{(x,x+dx)}(\xi_1) \right] :=   {d}P(x) = p(x \mid u) \, dx  
\] 
where $u$ is the parent's position. In a general BRW, offspring form a point process with a random number of offspring $\mathcal{N}$.   
By Wald's lemma \cite{Wald1945}, the corresponding intensity measure of the point process becomes
\[
       d Z(x \mid u) := \bold E \Big[ \sum_{j=1}^\mathcal{N}  \mathbf{1}_{(x,x+dx)}(\xi) \Big] 
     = \bold E [\mathcal{N}] \cdot  \bold E \left[ \mathbf{1}_{(x,x+dx)}(\xi) \right]  
     = m(u) \cdot  {d}P(x \mid u)  
\] 
or equivalently 
\[
     Z(X < x \mid u) := \bold E \Big[ \sum_{j=1}^\mathcal{N}   \mathbf{1}_{(-\infty,x)}(\xi) \Big] 
                      = \bold E [\mathcal{N}] \cdot  \bold E \left[ \mathbf{1}_{(-\infty,x)}(\xi) \right]  
     = m(u) \cdot P(X < x \mid u) 
\] 
where $   dP(x | u) = p(x | u)dx$ is the displacement law of a single offspring. The second equality above holds due to the independence of $\mathcal{N}$ and the offspring displacements. Here $m(u) = \bold E[\mathcal{N}]$ is the expected offspring (branching factor) by a parent at position $u$. In the above definition, the dropped subscript ``1'' should be replaced by the generation index $i$, which we omit for notational simplicity; in a general BRW, branching factor and the law of the offspring point process may vary across generations.

This measure $  d Z(x \mid u)$ at each generation uniquely determines the law of the BRW. Specifically, the probability law of the spine (a single lineage of descendants) is 
\begin{equation}\label{law of spine} 
             \frac { \displaystyle 
                     \prod_{i=1}^{n} m_i(u_{i-1}) p_i(x_i \mid u_{i-1}) \,\renewcommand{\thefootnote}{$\dagger$} \footnotemark 
                   }
                   { \displaystyle 
                     \int \prod_{i=1}^{n} m_i(u_{i-1}) p_i(x_i \mid u_{i-1}) \, dx_1 dx_2 \cdots dx_n
                   } 
                   \, dx_1 dx_2 \cdots dx_n 
\end{equation}
which can be interpreted as the proportion of spines $(x_1, x_2, ..., x_n)$ in the whole forest. 
\begingroup
\renewcommand{\thefootnote}{$\dagger$}
  \footnotetext{$\displaystyle u_{i-1} = x_0 + x_1 + \cdots + x_{i-1}$. 
               }
  \addtocounter{footnote}{-1}%
\endgroup                      

In view of this, a BRW can be defined by a sequence of pairs
\[
      \Big( m_i(u),\, p_i(x \mid u) \Big)_{i=1, 2, ...  } 
\]
where $m_i(u)$ is the expectation of offspring (branching factor) of a parent at position $u$, and $p_i(x \mid u)$ is pdf of  offspring's displacement, the probability (or proportion) density function, or, in discrete cases, mass probability. This definition covers a more general cases than the traditional setting; for instance, the parents' lifetimes need not be specified; whether they die or live, all such information is encoded in $m(u)$.

As mentioned earlier in this section, we are concerned with BRWs in which birth-rate is independent of birth position; that is, $m(u)$ does not depend on $x$. In this case, the law of spine (\ref{law of spine}) reduces to
\begin{equation}\label{law of RW}
           p_{1}(x_1) p_{2}(x_2) \cdots p_n(x_n)\, dx_1 \, dx_2 \cdots dx_n
\end{equation}
which is precisely the probability measure (law) for the random walk $(X_1, X_2, ..., X_n)$ where 
$$p_{i}(x) \equiv p(x| u_{i-1}) $$ 
This observation turns questions about BRW into questions about a random walk (without branching).

To aid intuition, we illustrate the new BRW framework in the discrete case. Let $M_0$ be the number of initial ancestors (can be infinity), and $M_i(u)$ the total population of $ith$ generation at position $u$, in the BRW forest.
The population, at $i+1$ generation, produced by those $M_i(u)$ particles is calculated by
\[
    M_i(u) \sum_x m_{i}(u)p_{i+1}(x)
\]
where
\[
    M_i(u)\cdot m_{i}(u)p_{i+1}(x)
\]
is population at $u+x$, $m_{i}(u)$ is the mean of birth-rate of a (parent) particle at $u$ in $ith$ generation, and $\sum_x p_{i+1}(x) = 1$. Note, $p_{i+1}(x)$ is dependent of $u$ which we drop off for notational simplicity; given $u$, $p_{i+1}(x)$ is viewed as a function of $x$ only, though. With
\[
    M_0 \cdot m_1(u_1)p_2(x)   ~~~~~ \mbox {($u_1=X_1$)}
\]
being the population (of 2nd generation) at $u_1 + x$ produced by 1st generation, we have by induction the total population at generation $n$ \big(note, $p_1(u_1) = 1$ \big)
\[
   \sum_{x_1}\sum_{x_2}\cdots\sum_{x_n}
        M_0 p_{1}(x_1) m_{1}(u_1)p_{2}(x_2)  \cdots m_{n-1}(u_{n-1})p_n(x_n)
\]
and the law of the spine $(x_1, x_2, ... , x_n)$

\begin{equation}\label{law of BRW}
 \frac{
          p_{1}(x_1) m_{1}(u_1)p_{2}(x_2)  \cdots m_{n-1}(u_{n-1})p_n(x_n)
        }
        {
            \sum_{x_1}\sum_{x_2}\cdots\sum_{x_n}
                  p_{1}(x_1) m_{1}(u_1)p_{2}(x_2)  \cdots m_{n-1}(u_{n-1})p_n(x_n)
        }
\end{equation}
which is the proportion of spines $(x_1, x_2, ..., x_n)$ in the whole forest. If the birth-rate is independent of birth position, $m(u)$ does not depend on $x$. In this case, the law of the spine (\ref{law of BRW}) reduces to
\begin{equation}\label{law of RW}
           p_{1}(x_1) p_{2}(x_2)\cdots p_n(x_n)
\end{equation}
which is the probability measure of the random walk $(X_1, X_2, ..., X_n)$ where $$ p_{i}(x) \equiv \Pr(X_i=x \ | \ u_{i-1}) $$ 
   
{\textbf {Remarks.}}

\textbullet\ As far as the stochastic behavior at generation $n$ is concerned, random walk can be viewed as a special BRW (with branching factor = 1), or random walk is a special BRW where birth-rate is independent of birth-place.

\textbullet\ A BRW of spatial homogeneity in branching (i.e. $m_n(u)=m_n$) can be studied as a random walk process (without branching). In particular, 
\[ 
       \int_{\left|\, S_n - \bold E(S_n) \, \right| \ge \lambda}
               \frac {   
                       \prod_{i=1}^{n} m_i(u_{i-1}) p_i(x_i \mid u_{i-1}) 
                     }
                     {   
                       \int \prod_{i=1}^{n} m_i(u_{i-1}) p_i(x_i \mid u_{i-1}) \, dx_1 dx_2 \cdots dx_n
                     } 
                     \, dx_1 dx_2 \cdots dx_n 
\]
reduces to a large deviation problem for the underlying random walk, namely
\[       
     \int_{\left|\, S_n - \bold E(S_n) \, \right| \ge \lambda}
                                \prod_{i=1}^{n}  p_i(x_i \mid u_{i-1})  \, dx_1 dx_2 \cdots dx_n 
                             = \Pr \left( | S_n - \bold E(S_n)| \ge  \lambda  \right)
\]

In the following, we denote the random walk by $(p_{i})_{i=1, 2, ...n}$ \big(note this random walk is not necessarily the underlying random walk of the original BRW, though under the assumption of independence between siblings' motion, the underlying random walk is $(p_{i})$\big). The following theorem presents a Chernoff bound for branching random walks; the classical i.i.d. model of Harris (1960) arises as a special case.
\begin{thm}\label{BRW_concentration} For BRW $(m_i, p_i)_{i=1, 2, ..., n}$, if 
\\
\indent(a) $m_i$ is position-independent, and
\\
\indent(b) Chernoff bound holds for the random walk $(p_i)$, then 
\[
       \Pr  \big( |Q_n(\alpha) - n\mu | \ge \lambda \big) \le 2e^{-\frac{1}{2}c\lambda^2/n}
\] 
where $n\mu$ is the expected position\footnote{``$n\mu$" is used for ``mean'' in order to be consistent with the literature of BRW, not indicating a linear relationship with $n$} for an individual in $nth$ generation, $Q_n(\alpha)$ is $\alpha$ quantile and 
\[ 
       \alpha \in \left[ e^{-\frac{1}{2}c\lambda^2/n}, 1-e^{-\frac{1}{2}c\lambda^2/n} \right]
\] 
Throughout, generation index $n$ will be omitted when no confusion can arise.
\end{thm}
The theorem conveys the same spirit as the Chernoff inequalities for random walks, stating, roughly, that with high probability (with exception probability of order $e^{-O(1)\lambda^2/n}$), at least a proportion $1 - e^{-O(1)\lambda^2/n}$ of the leaves of a BRW tree are concentrated within a distance of order $\sqrt{n}$ from the expectation, in the regime $\lambda \gtrsim \sqrt{n}$.

\begin{proof}
Let $z^{(n)}_1, z^{(n)}_2, ...$ be an enumeration of the positions of the particles (leaves) in the $nth$ generation and $Z^{(n)}$ be its population; i.e. $Z^{(n)} = |\{z^{(n)}_1, z^{(n)}_2, ...\}|$. There should be an index variable $\tau$ for trees in the above notations which we omit. Consider
\[
    \frac{\sum_i  \mathbf{1}_{\{z_i \le t\}} }
           {Z^{(n)}}
\]
namely, the cumulative distribution function of the data set $\{z^{(n)}_1, z^{(n)}_2, ...\}$.
The $1-\alpha$ quantile is
\[  Q_n(1-\alpha) = \inf \{t: 
       \frac{\sum_i  \mathbf{1}_{\{z_i \le t\}} }
            {Z^{(n)}} 
            \ge 1-\alpha  \}
\]

Let $\tau$ be a tree of $Z$ leaves and $p_{(\tau)}(\lambda)$ be the proportion of positions which $\ge n\mu + \lambda $ (in the $n$th generation), i.e.
\[
    p_{(\tau)}(\lambda) = \frac {1} {Z} \sum_{i=1}^Z \mathbf{1}_{ \{  z_i(\tau) \ge n\mu + \lambda \}}
\]
By definitions,
\begin{equation}\label{definitions}
      p_{(\tau)}(\lambda) \ge \alpha  \iff  Q_n(1-\alpha) \ge  n\mu + \lambda     
\end{equation}
We have 
\[
    \bold{E}\left[ p_{(\tau)}(\lambda) \right] = \Pr(S_n - n\mu \ge  \lambda) 
\] 
because
\begin{align}\label{Ep_tau} 
   \mathbf{E} \left[ p_{(\tau)}(\lambda) \right] = \frac {\sum_{i=1}^Z \bold E(\mathbf{1}_{ \{  z_i - n\mu \ge \lambda \}}
                                                                              )
                                                         }
                                                         {Z}
                          = \frac {\sum_1^Z \Pr(S_n - n\mu \ge  \lambda) }{Z}
                          = \Pr(S_n - n\mu \ge  \lambda)
\end{align}
where $\tau$ is a tree of population $Z$ (of generation $n$). In other words, the expectation of $p_{(\tau)}(\lambda)$ over all the trees of the same population is 
\[
   \bold{E}\left[ p_{(\tau)}(\lambda) \right]  =  \Pr(S_n - n\mu \ge \lambda)      
\]
The reason $Z$ does not appear explicitly in this expectation is that, for a randomly selected leaf $z_i$ from a tree, 
\[ \bold{E} \left[\mathbf{1}_{ \{  z_i - n\mu \ge \lambda \}} \right]
\] 
is independent of the size of the tree because, by the hypothesis $(a)$, branching is independent of walking. Branching is also independent from leaf indexing so that $E(\mathbf{1}_{ \{  z_i - n\mu \ge \lambda \}})$ is the same for any $i$. Below we will use $k$ instead of capital $Z$ for readability in summation.


On the other hand, denoting the number of trees of size $k$ by $n_k$, we have
\begin{align}\label{sampling of E}
     \frac{1}{k}  {\sum_{i=1}^k \bold{E} \left[ \mathbf{1}_{ \{  z_i - n\mu \ge \lambda \}} \right] }    =
          \frac{1}{k} \sum_{i=1}^k \frac{1}{n_k}\sum_{\tau_k} \mathbf{1}_{ \{  z_i(\tau_k) - n\mu \ge \lambda \}}
\end{align}
Note that $n_k$ is very large, since we have infinitely many initial ancestors. Suppose there are $N_t$ trees in the forest, $n_1$ trees of population 1, $n_2$ trees of population 2, ... and so on.
In view of above equations, we have
\begin{align*}
      \Pr(S_n - n\mu \ge  \lambda)  &= \sum_{k\ge 1}    \frac{n_k}{N_t}    \frac{1}{k} \sum_{i=1}^k  
      \frac{\sum_{\tau_k} \mathbf{1}_{ \{  z_i(\tau_k) - n\mu \ge \lambda \}} }
                         {
                          n_k
                         }
\mbox {~~~ (\ref{Ep_tau}) and then (\ref{sampling of E})
                        }
                        \\
               &=
 \frac{1}{N_t} \sum_{k\ge 1}   \sum_{\tau_k}  \frac{\sum_{i=1}^k \mathbf{1}_{ \{  z_i(\tau_k) - n\mu \ge \lambda \}} }
                                                     {k}  \\
               &=
\frac{1}{N_t}    \sum_{\tau}  \frac{\sum_{i=1}^k \mathbf{1}_{ \{  z_i(\tau) - n\mu \ge \lambda \}} }
                                                     {k}               \\
               &=
\frac{1}{N_t}    \sum_{\tau}  p_{(\tau)}(\lambda)   \\
               & \ge \frac{1}{N_t}    \sum_{\tau}  p_{(\tau)}(\lambda) \mathbf{1}_{\{p_{(\tau)}(\lambda) \ge \alpha\}}
                                                 \\
               & \ge \frac{1}{N_t}    \sum_{p_{(\tau)}(\lambda) \ge \alpha }  \alpha  
                                                   =  \alpha \cdot  \Pr\Big(p_{(\tau)}(\lambda) \ge \alpha \Big)  
                 \\
               & =  \alpha \cdot  \Pr\Big(Q_n(1-\alpha) \ge n\mu + \lambda \Big) \quad  \text{by (\ref{definitions})}
\end{align*}
By the well known Chernoff bound for the left-hand side of the above, it follows, 
\[
       \alpha \cdot  \Pr \Big(Q_n(1-\alpha) \ge n\mu + \lambda  \Big)   \le  e^{-c\lambda^2/n}
\]
where $c > 0$ is a constant that may depend on $n$. 
Choosing $\alpha = e^{-\frac{1}{2}c\lambda^2/n}$, we have
\begin{equation}
                       \Pr\Big(Q_n(1-e^{-\frac{1}{2}c\lambda^2/n}) - n\mu \ge \lambda \Big) \le e^{-\frac{1}{2}c\lambda^2/n}
\end{equation}
Since $Q_n(\alpha)$ is monotonically non-decreasing in $\alpha$, it follows that for all 
$\alpha \le 1- e^{-\frac{1}{2}c\lambda^2/n}$,%
\begingroup 
\renewcommand{\thefootnote}{$\S$}
\footnote{Here the symbol $\alpha$ is reused to denote a generic parameter, rather than the specific value chosen above.}
     \addtocounter{footnote}{-1}     
\endgroup%
\begin{equation}
                       \Pr(Q_n(\alpha) - n\mu \ge \lambda   ) \le e^{-\frac{1}{2}c\lambda^2/n} \, ;
\end{equation}

Similarly, for $\alpha \ge e^{-\frac{1}{2}c\lambda^2/n}$
\begin{equation}
                       \Pr \Big( Q_n(\alpha) - n\mu \le  - \lambda  \Big) \le e^{-\frac{1}{2}c\lambda^2/n}
\end{equation}
The claim follows.

\end{proof}

{\textbf {Remarks.}}

\textbullet\ The Chernoff bound holds for both supercritical and subcritical BRW, because it is independent of the branching factor $m_i$. 

\textbullet\ The Chernoff bound presented in this paper is a special case of the more general concentration inequalities (\cite{Liu Con Ineq}), which holds in a setting where (a) $m_i(u)$ is position-dependent (and thus path-dependent), (b) neither independence nor martingale conditions are assumed for the underlying random walk processes, and (c) the boundedness assumption on the step increments $X_i$ is relaxed. 

\textbullet\ If displacements of siblings are independent of each other, with step size $X_i$ then, $p_i$ and the law of $X_i$ are identical (note, generally they are not equal). In other words, if $S_n = X_1 + \cdots + X_n$ is a random walk with $(p_i)$ as the probability density of the increment, then the distribution of $S_n$ has the same "shape" as the forest of BRW $(m_i, p_i)$. Bear in mind, probability space of random walk is interpreted as the forest of BRW of $(1, p_i)$.

\textbullet\ In CSPs (constraints satisfaction problem) such as $K$-SAT, the enumeration of the problem instances can be formulated as a branching random walk \cite{Liu K-SAT/q-COL}. In this BRW, the forest has only one tree because every tree is the same; $m_i$ is not random given it birthplace, and in addition branching factor is large (say $(2n)^k$). Because the branching structure depends on the position, the BRW can not be reduced to random walk;  (\ref{law of BRW}) can not be reduced to (\ref{law of RW}). The concentration inequalities for this BRW are developed in a separate paper \cite {Liu Con Ineq}.


\end{document}